\pdfoutput=1

\documentclass{amsart}
\usepackage[centertags]{amsmath}
\usepackage{amsfonts}
\usepackage{amssymb}
\usepackage{amsthm}
\usepackage{subfigure}
\usepackage{fullpage }
\usepackage{color}
\usepackage{graphicx}

\newtheorem{theorem}{Theorem}
\newtheorem{proposition}[theorem]{Proposition}
\newtheorem{lemma}[theorem]{Lemma}

\theoremstyle{definition}
\newtheorem{definition}[theorem]{Definition}

\theoremstyle{remark}

\begin{document}

\renewcommand{\labelenumi}{(\roman{enumi})}

\def\printname#1{
	
		\smash{\makebox[0pt]{\hspace{-0.5in}
			\raisebox{8pt}{\tt\tiny #1}}}
}
\def\lbl#1{\label{#1}
}

\newcommand{\fD}{\mathfrak{D}}
\newcommand{\fG}{\mathfrak{G}}
\newcommand{\fL}{\mathfrak{L}}
\newcommand{\fK}{\mathfrak{K}}
\newcommand{\fs}{\mathfrak{s}}

\newcommand{\Z}{\mathbb{Z}}
\newcommand{\D}{\mathbb{D}}
\newcommand{\G}{\mathcal{G}}
\newcommand{\M}{\mathcal{M}}
\newcommand{\R}{\mathcal{R}}
\newcommand{\p}{\mathcal{P}}
\newcommand{\A}{\mathcal{A}}
\newcommand{\E}{\mathcal{E}}
\newcommand{\V}{\mathcal{V}}
\newcommand{\s}{\mathcal{S}}
\newcommand{\F}{\mathcal{F}}

\newcommand{\BR}{Bollob\'{a}s-Riordan }


\title[ Unsigned state models for the Jones polynomial]{ Unsigned state models for the Jones polynomial}

\author[I.~Moffatt]{Iain Moffatt$^*$}

\begin{abstract}
It is well a known and fundamental result that the Jones polynomial can be expressed as Potts and vertex partition functions of signed plane graphs. Here we consider constructions of the Jones polynomial as state models of unsigned graphs and show that the Jones polynomial of any link can be expressed as a  vertex model of an unsigned embedded graph.

In the process of deriving this result, we show that for every diagram of a link in $S^3$ there exists a diagram of an alternating link in a thickened surface (and an alternating virtual link) with the same Kauffman bracket. We also recover two recent results in the literature relating the Jones and \BR polynomials and show they arise from  two different interpretations of the same embedded graph.  
\end{abstract}

\thanks{
${\hspace{-1ex}}^*$
Department of Combinatorics and Optimization,  University of
Waterloo, Waterloo, Ontario, Canada; \\
${\hspace{.35cm}}$ {\em Current address:} Department of Mathematics and Statistics,  University of South Alabama, Mobile, AL 36688, USA;\\
${\hspace{.35cm}}$ \texttt{imoffatt@jaguar1.usouthal.edu}
\\
This version: April 28, 2009. First version: October 20, 2007.
}

\maketitle
\section{Introduction and motivation}\lbl{sec:int}
Connections between knot theory and statistical mechanics were first noted over twenty years ago by Jones \cite{Jo2}, and were made concrete by Kauffman \cite{Ka,Ka2} shortly after. Such relations between physics and knot theory were soon found for other knot invariants by Jones  \cite{Jo} and  Turaev \cite{Tu}, and have since been explored  by many others. We refer the reader to either  Jones' paper \cite{Jo} or Wu's survey article \cite{Wu1} for an overview  of connections between statistical mechanical models and knot invariants.
 
 In this paper we are interested in expressions of the Jones polynomial as a Potts (or  spin) model and as a vertex (or ice-type) model and also the relations of these models to graph polynomials. These two models can be roughly described as follows. 
 Given a oriented link diagram $D$, there are two standard ways to   construct signed graphs from $D$.  Each of these ways gives rise to a statistical mechanical model for the Jones polynomial. 
 
  Firstly, we can regard the diagram $D$ as a four-valent plane graph itself and add a sign $+$ or $-$ to each vertex of the graph according to the (oriented) sign of the corresponding crossing. This gives a vertex weighted graph $G_D$. The Jones polynomial  of $D$ can be given as a  sum over all edge-colourings of $G_D$  where the summands in this state sum arise from the graph colouring locally at each vertex. Such a construction is known as a vertex model.
 
 On the other hand, we can consider the Tait graph of $D$ (this is also known as the medial graph in the literature, but here we will use the term ``medial graph'' in its graph theoretical sense). The Tait  graph $M_D$ of $D$  is an edge-signed plane graph obtained from $D$ by considering the checker-board colouring of $D$; assigning vertices to the black regions and edges where black regions meet (see Subsection~\ref{sec:kauffman} for details). The Jones polynomial can then be obtained as a state sum over all spanning subgraphs of $M_D$. The summands in this state sum  are determined by the signs of the edges in the spanning subgraph. A statistical mechanical construction of this type is known as a Potts model.   
 
Relating to the above information, we observe the following facts.
  The vertex model construction of the Jones polynomial from \cite{Tu} can be recovered from the Potts model construction of the Jones polynomial in \cite{Ka3} using the idea of arrow coverings from  \cite{Ba}.  The Kauffman bracket arises as an intermediate step of such a recovery of a vertex model from the Potts model.  Also  Thistlethwaite's theorem \cite{Th, Ka3} connecting the Jones and Tutte polynomials can be recovered from the Potts model construction of the Jones polynomial. 
 These facts will motivate the results presented in this paper.

Having obtained one collection of motivational facts, we will now set about  acquiring  another collection. Comparing these two collections of facts will motivate the results presented here.
 
  \medskip
 
 Recently there has been interest in connections between the \BR polynomial and knot invariants (for example in \cite{CP,CV,Da,LM,HM,Mo}).
 The \BR polynomial  of a ribbon graph (\cite{BR1, BR}) is a  generalization of  the Tutte polynomial to ribbon graphs or embedded graphs. The \BR polynomial depends  upon the topology of the ribbon graph in an essential way and the Tutte polynomial can be recovered as a specialization of the \BR polynomial.  Here we are particularly  interested in a result from
 \cite{Da}, where  Dasbach,  Futer,  Kalfagianni,  Lin and  Stoltzfus showed how to construct a ribbon graph $\D$ from an arbitrary link diagram $D$, and showed how the Kauffman bracket $\langle D \rangle$ can be obtained as an evaluation of the \BR polynomial  of the ribbon graph $\D$. We emphasise the fact that this is a calculation of the Jones polynomial of a (not necessarily alternating) link as a  polynomial of  an unsigned (but not necessarily planar) graph. 
When the link is alternating, the ribbon graph  constructed in \cite{Da} is exactly the Tait graph $M_D$ of the link diagram. The relation between the \BR polynomial of $\D$ and the Jones polynomial of $D$ specialises to a relation between the Tutte polynomial of $M_D$ and the Jones polynomial of $D$. Therefore Dasbach,  Futer,  Kalfagianni,  Lin and  Stoltzfus' result is a generalization of Thistlethwaite's theorem relating the Tutte and Jones polynomials.

This gives rise to our second collection of motivational facts: the Jones polynomial of a link can be obtained as an evaluation of the \BR  polynomial of an unsigned ribbon graph, and this connection is a generalization of Thistlethwaite's  connection between the Jones polynomial and  the Tutte polynomial.

\medskip

Our two  collections of facts  motivate some questions. The first question is straightforward. Given  the close connection between Tutte polynomial and the Potts models, and the fact that the \BR polynomial generalises the Tutte polynomial, how do the results of \cite{Da} relate to the Potts model for the Jones polynomial? Or equivalently, how do unsigned state models for the Jones polynomial relate to the signed state models for the Jones polynomial? Secondly, if we regard  the results from \cite{Da} as, in some sense, a Potts model construction of the Jones polynomial using unsigned graphs (this is reasonable given the close connections between the Potts partition function, the Tutte polynomial and the \BR polynomial), can we use this to construct a vertex model for the Jones polynomial which uses unsigned graphs?

Here we address these questions.  Starting with the Potts model for the Kauffman bracket we show that  the Kauffman bracket (and therefore the Jones polynomial) of a link can be obtained as a vertex model of an unsigned embedded graph.  
During  this process we  recover results \cite{Da, CP},  which relate the \BR polynomial and the Jones polynomial of links and links in thickened surfaces, and show that they arise from two different interpretations of the same embedded graph. We also show for every link in $S^3$  there is an alternating  link in a thickened surface (or 
a  virtual alternating link) which has the same Kauffman bracket.

The approach we use in this paper  is motivated by  the method  from \cite{Ba} used to show that the Potts model can be  expressed as a vertex model, and the construction of the ``all-A ribbon graph'' from \cite{Da}. 

\smallskip

I would like to thank the referee for his or her helpful comments.

\section{Preliminaries}\lbl{sec:pre}
\subsection{Graphs and ribbon graphs}

Let $G=(V,E)$ be a  graph (possibly with multiple edges and loops) with vertex set $V$ and egde set $E$. If $A\subset E$, then the subgraph $(V,A)$ of $G$ is called a spanning subgraph. 
We let  $k(A)$  denote the number of connected components, $r(A):=V-k(A)$ denote the rank, and $n(A)= A-r(A)$ denote the nullity 
of the spanning subgraph $(V,A)$.

A graph $G=(V,E)$  is said to be edge (respectively vertex) $q$-coloured if there is a map from $E$ (respectively $V$) to the $q$ element set $\{1, 2, \ldots , q\}$. We say that a graph is {\em signed} if  
there is a map from its vertex or edge set to the two element  set $\{ +, -\}$.
 If $G=(V,E)$ is signed then we let $E_{\pm}$ denote the set of $\pm$ coloured edges of $G$. Note that  $E_+$ and $E_-$ partition $E$.

\medskip

A {\em ribbon graph} (or {\em combinatorial map} or {\em fatgraph}) is a graph
(possibly with multiple edges and loops) with a fixed cyclic ordering
of the incident half-edges at each of its vertices.

A ribbon graph can be realized as an orientable surface by
``thickening'' its edges and vertices. The following is clearly an
equivalent definition of ribbon graphs. A ribbon graph $\F=(V(\F),E(\F))$ is
an orientable surface with boundary represented as the union of
$|V(\F)|$ closed {\em discs} and $|E(\F)|$ {\em ribbons}, $I \times I$, such that
\begin{enumerate}
\item the discs and ribbons intersect in disjoint line segments
$\{0,1\} \times I$;
\item each such line segment lies on the boundary of precisely one
disc and precisely one ribbon;
\item every ribbon contains exactly two such line segments.
\end{enumerate}
 When we want to emphasise the surface nature of a ribbon graph, we  refer to edges as ribbons
and vertices as discs.

Ribbon graphs can also be regarded as graphs embedded in an orientable surface as follows. The natural embedding of a ribbon graph in its ``thickening'' or ``surface realization''  provides an embedding of a graph in an orientable punctured surface. We can then obtain an embedded graph in a closed orientable surface $\Sigma$ by capping  the punctures using discs.

We thus have three perspectives of a ribbon graph and we can move freely between these three notions at our convenience. 
An example of the three realizations of a ribbon graph is shown in Figure~\ref{fig:rib}.
\begin{figure}
\includegraphics[width=3cm]{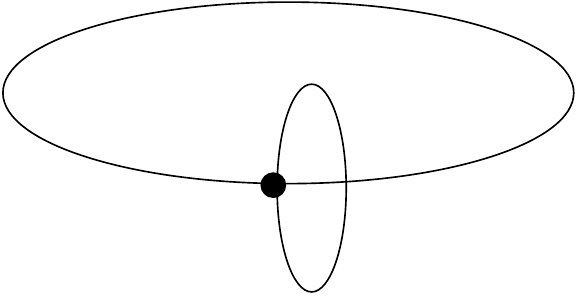}
 \hspace{2cm} 
 \includegraphics[width=2cm]{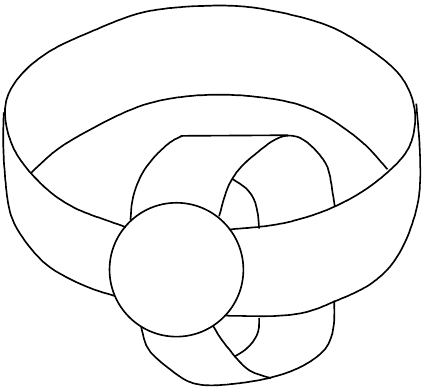}
 \hspace{2cm} 
 \includegraphics[width=3cm]{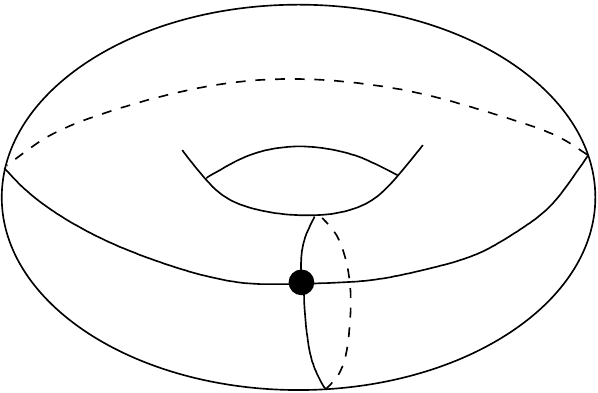}
\caption{Three realizations of a ribbon graph}
\lbl{fig:rib}
\end{figure}

If $\F= (V(\F) , E(\F))$ is a ribbon graph then $\A \subset E(\F)$ determines a spanning ribbon sub-graph $(V(\F) , \A)$ of $\F$. Since a ribbon graph consists of a graph with some additional structure, graph theoretical concepts, such as colourings, connectivity, rank and nullity, also apply to ribbon graphs. 
In addition to these graph theoretical concepts, by regarding a ribbon graph as a surface, we let
$\partial(\A)$  be the number of boundary components of   $(V(\F) , \A)$. 
We note that if   $\F=(V(\F),E(\F))$ is a ribbon graph viewed as a graph embedded in a surface $\Sigma$, we regard its spanning ribbon sub-graphs  as graphs which are also embedded in the surface $\Sigma$.

\subsection{Statistical mechanical models}
We are interested in  vertex models and Potts models - two of the fundamental models in statistical mechanics. These two models are defined as follows.

Let $G=(V,E)$ be an edge q-coloured graph. Let $\sigma$ denote the edge q-colouring of $G$.
For each vertex $v\in V$ we assign a {\em weight function} $w_v(\sigma)$ which takes its value in some ring. The value of  $w_v(\sigma)$ is determined by the edge q-colourings of the edges incident with $v$. (The $w_v$ are known as the Boltzmann factors.)
We define a {\em q-state vertex model} to be a graph $G$ equipped with a weight function for each vertex.  The {\em partition function} for a q-state vertex model is the state sum
\begin{equation}\label{eq:ver}
Z_{vertex}(G) = \sum_{\sigma: E \rightarrow \{0, 1, \ldots , q-1\}}\prod_{v\in V}w_v(\sigma).
\end{equation}
The sum is over all edge q-colourings $\sigma$ of $G$.

\medskip

The Potts model and its partition function also arise  from vertex weighted graphs, however the weight functions in this model  depend upon the colouring of the vertices incident with each edge of $G$, rather than the colouring of the  edges incident with each vertex.
 
Let $G=(V,E)$ be an vertex q-coloured graph.
We define a {\em q-state Potts model} to be a graph $G$ such that each edge $e$ is equipped with a weight $w_e$ .  The {\em Potts partition function} 
is given by\begin{equation}\label{eq:pot}
Z_{Potts}(G) = \sum_{\sigma: v \rightarrow \{0, 1, \ldots , q-1\}}\prod_{e\in E}
\left[1+ w_e\delta \left(\sigma({v_1}) , \sigma({v_2})  \right) \right],
\end{equation}
where $v_1$ and $v_2$ are the vertices of the edge $e$, and $\delta(a,b)$ is the Kronecker delta.

The Potts model also has the following realization, which we will find  convenient. 
\begin{theorem}[\cite{FK}] \lbl{th:forkas} 
Let $G=(V,E)$ be a graph. Then for each $q\geq 1$, we have
\[ 
Z_{Potts}(G)=\sum_{A\subset E }  q^{k(A)}\prod_{e\in A}w_e ,
\]
\end{theorem}
(See \cite{Ba} or \cite{So} for the short proof of the theorem.)

The state sum in the theorem will be important to us, so we shall set 
\[  Z(G) =\sum_{A\subset E }  q^{k(A)}\prod_{e\in A}w_e.  \]
 The state sum  $Z(G)$  is  known in physics as the Fortuin-Kasteleyn representation of the Potts model, and   in graph theory as the  mulitvariate Tutte polynomial \cite{So}. 

For certain classes of graph, the Potts and vertex models are equivalent (see \cite{Ba,PeWu}). As we have mentioned previously, we are  interested in the approach for expressing the Potts model of a plane graph as a vertex model using arrow coverings as in \cite{Ba}. Our ``unsigning'' of the Potts model arises from this approach.

\section{From signed graphs to ribbon graphs}\lbl{sec:rewrite}
In this section we show that, provided that the edge-weights satisfy the identity $w_+=qw_-^{-1}$, the Potts partition function of a signed plane graph $G$ can be written as a state sum  of an unsigned ribbon graph $\R_G$, and as a state sum of an  unsigned graph $\M_{G}$  embedded in an orientable  surface.

\subsection{The unsigned ribbon graph $\R_G$}
Let $G=(V,E)$ be an edge-signed plane graph. 
We construct an unsigned ribbon graph  $\R_G$ from $G$ by  replacing each signed edge of $G$ with a ribbon-arc configuration as indicated: 
%
\[ 
\includegraphics[height=1.5cm]{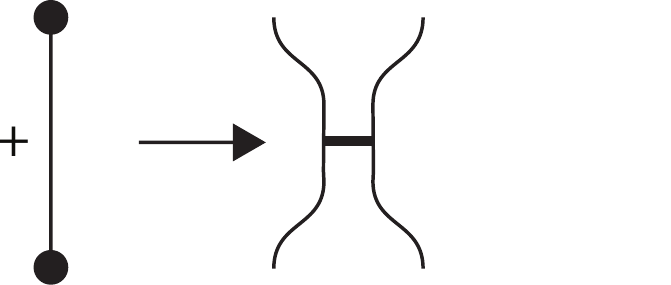}
 \hspace{2cm} 
 \includegraphics[height=1.5cm]{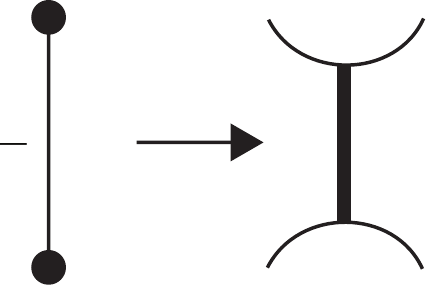}\]
and then by joining up these configurations according to the cyclic order of the incident edges of $G$. (In the figure the ribbons are shown as thick black lines.)
What we end up with is a collection of planar cycles  joined by lines. We may regard such a collection as a ribbon graph, by regarding each cycle as the boundary of a disc.  Form $\R_G= (V(\R_G) , E(\R_G))$, by taking $V(\R_G)$ to be the set of discs determined by the planar cycles, and whenever two  planar cycles are joined by a line, add a ribbon in the position specified, to obtain $E(\R_G)$.  
A example is shown in Figure~\ref{fig:unsignexamp}.
 We call the ribbon graph $\R_G$ thus obtained, the {\em unsigning of $G$}. 
\begin{figure}
\includegraphics[height=2cm]{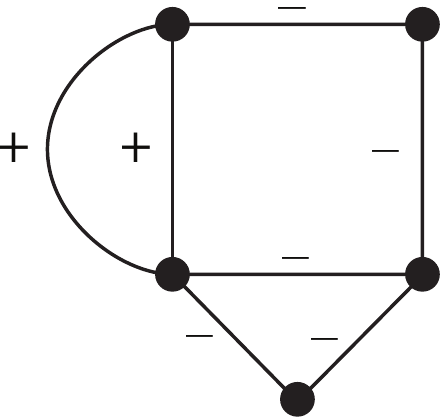}
 \hspace{2cm} 
 \includegraphics[height=2.5cm]{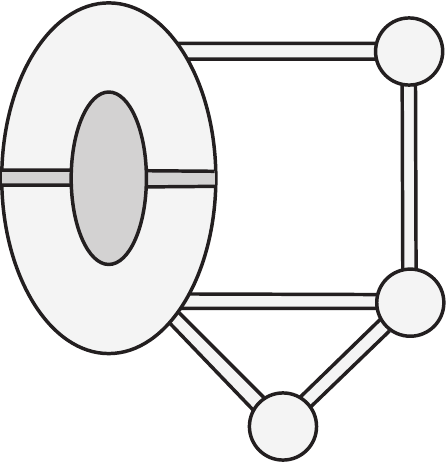}
\caption{An example of unsigning a plane graph}
\lbl{fig:unsignexamp}
\end{figure}

We note that if $G$ is the Tait graph of a link diagram then $\R_G$ is exactly the all-A ribbon graph from \cite{Da}. We will return to this point  in Section~\ref{sec:kauffman}.

We can write the Potts partition function for $G$ as a state sum of $\R_G$: 
\begin{proposition}\lbl{pr:unsig}
Let $G=(V,E)$ be an edge-signed plane graph and  $\R_G = (V(\R_G), E(\R_G))$ be its unsigning constructed as above. Then if    $x_+=(x_-)^{-1}$, where $x_{\pm}=q^{-1/2}w_{\pm}$, we have
\[   Z_{Potts}(G) = q^{V/2} x_+^{|E_+|} \sum_{ \A \subset E(\R_G)  } q^{\partial (A) /2} x_-^{|A|}  ,\]
where $|E_+|$ is the number of positively weighted edges in $E(G)$.
\end{proposition}

The remainder of this subsection is devoted to the proof of this proposition.

\medskip


We may regard the signed plane graph $G$ as a signed ribbon graph $\G=(V(\G),E(\G))$ by ``fattening'' the vertices  to discs and the edges into ribbons as indicated:
\[
\includegraphics[height=1.5cm]{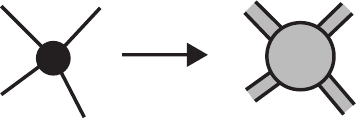}
.\]
The colouring of the ribbons of $\G$ is induced from the colouring of the edges of $G$.
We call the signed ribbon graph $\G$ the {\em fattening} of $G$.

By Theorem~\ref{th:forkas}, we  can write $Z_{Potts}(G)$ as the state sum
\[  Z(G) = \sum_{A\subset E }  q^{k(A)}\prod_{e\in A}w_{e}.  \] 
This state sum for $G$ can clearly be written in terms of the fattening $\G$ of $G$: 
\[   Z(G) =  \sum_{A \subset E(\G)}  q^{k(A)}\prod_{e\in A}w_{e}  .   \]

Let $(V,A)$ and $(V(\G) , \A)$ be spanning subgraphs of $G$ and $\G$ respectively.
Since $G$ and $\G$ are embedded in the plane, the two spanning subgraphs are too. By Euler's formula, we have 
$ |V|-|A|+f(A)=k(A)+1$, where $f(A)$ is the number of faces of the embedded spanning subgraph $(V,A)$.   Since $f(A)=\partial (\A)-k(\A)+1$, it  follows that $2k(\A)+|\A|-|V(\G)|-\partial(\A)=0$.
Applying this observation to the expression for $Z(G)$ above,  we obtain the following lemma.
\begin{lemma}\lbl{lem:pp1}
Let $G$ be a signed plane graph and $\G=(V(\G), E(\G))$ be its fattening. Then
\begin{equation} \label{eq:ac1}
  Z(G) = q^{|V(\G)|/2}  \sum_{A \subset E(\G)}  q^{\partial (A)/2}\prod_{e\in A}x_{e},\end{equation}
where $x_e = q^{-1/2}w_e$.
\end{lemma}

To prove Proposition~\ref{pr:unsig}, we need to relate the the state sum (\ref{eq:ac1}) to the state sum for $\R_G$ stated in  the proposition. This relation will follow as a straightforward  application of   the following lemma.
\begin{lemma}\lbl{lem:pp2}
 Given a signed plane graph $G=(V,E)$. Let $\G$ be its fattening and $\R_G$ its unsigning. Then if $y_+=(y_-)^{-1}$,
 \begin{equation}\label{eq:ac2}
 \sum_{A\subset E (\G)} a^{\partial (A)} \prod_{e\in A} y_e = 
 (y_+)^{|E_+|} \sum_{\A \subset E(\R_G)} a^{\partial (\A)} (y_-)^{|\A|},
 \end{equation}
 where $|E_+|$ is the number of positive edges of $G$ (or $\G$).
\end{lemma}
\begin{proof}
We construct a bijection between the spanning subgraphs of $\G$ and $\R_G$ with the property that corresponding states contribute the same terms to the expressions on the left and right of Equation~(\ref{eq:ac2}).

Notice first of all that there is a natural correspondence between the ribbons of $\G$ and the ribbons of $\R_G$. We use this correspondence throughout the proof. 

Given $A\subset E(\G)$ construct a spanning subribbon graph $\s_A$ of $\R_G$ as follows:
\begin{itemize}
\item if $e\in A$ and  \begin{itemize} \item if $e$ is of positive sign remove the corresponding ribbon from $\R_G$;
\item if $e$ is of negative sign leave in the corresponding ribbon from $\R_G$;
\end{itemize}
\item if $e\notin A$ and  \begin{itemize} \item if $e$ is of negative sign remove the corresponding ribbon from $\R_G$;
\item if $e$ is of positive sign leave in the corresponding ribbon from $\R_G$.
\end{itemize}
\end{itemize}

First observe that  if $n_{\pm}$ denotes the number of $\pm$ edges in $A$, then $\s_A$ has $n_- +(|E_+| - n+)$ edges. If $y_+=(y_-)^{-1}$, we then have
\[   y_-^{|\s_A|} =  y_-^{n_- +(|E_+| - n+)}  = y_-^{|E_+|}  y_-^{n_-} y_-^{ - n+} = y_+^{-|E_+|}  y_-^{n_-} y_+^{  n+}   =  y_+^{-|E_+|} \prod_{e\in A} y_e .\]

Secondly, we need to show that $\partial(A) = \partial (\s_A)$.
To show this we begin by proving  that  the boundary components of $\G$ are in correspondence with the boundary components of $\s_{E(\G)}$. 
To prove this correspondence between boundary components, we start by considering the arcs in the boundary components of $\G$ arising at an edge $e$, and the  arcs of the boundary components of $\s_{E(\G)}$ arising at the corresponding ribbon-arc configuration. From the figure below, we see that locally at the edge $e$ of $\G$, the boundary components define arcs from $a_e$ to $b_e$ and from $c_e$ to $d_e$, as shown. Also, the boundary components at the corresponding edge of   $\s_{E(\G)}$ will also define arcs from $a_e$ to $b_e$ and from $c_e$ to $d_e$, as shown.

\begin{center}
\begin{tabular}{p{3cm} c p{3cm} c p{3cm} cp{3cm} } 
\includegraphics[height=2.5cm]{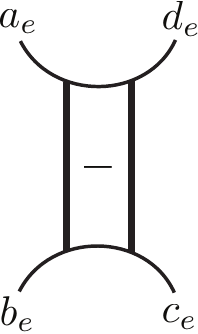}  &\hspace{-17mm}\raisebox{12mm}{\includegraphics[width=1cm]{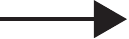}  }& \includegraphics[height=2.5cm]{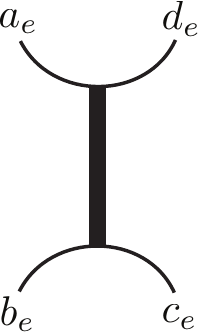}   & & \includegraphics[height=2.5cm]{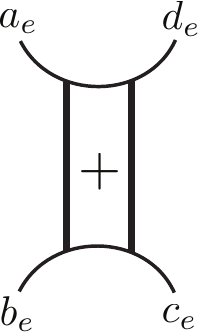}  & \hspace{-17mm}\raisebox{12mm}{\includegraphics[width=1cm]{figs/arrow}  }&\includegraphics[height=2.5cm]{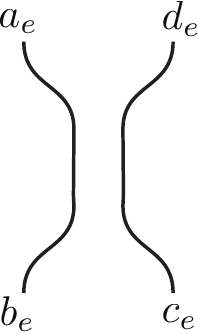}  
\\
edge $e$ in $\G$ & &corresponding configuration in $\s_{E(\G)}$ & \hspace{1cm} & edge $e$ in $\G$ && corresponding configuration in $\s_{E(\G)}$ 
\end{tabular}
\end{center}

\noindent By the definition of $\R_G$, the points in the set $\{ a_e, b_e, c_e, d_e\, | \, e\in E(\G)\}$ in $\G$ are connected  to each other in in exactly the same way that the points in the set $\{ a_e, b_e, c_e, d_e\, | \, e\in E(\s_{E(\G)})\}$ in $\s_{E(\G)}$ are connected  to each other. This means that the local configurations of $\G$ and $\s_{E(\G)}$ at each edge $e$ shown in the figure above are connected in the same way. The correspondence between the boundary components of $\G$ and $\s_{E(\G)}$ then follows.

We will now consider the effect of the removal of ribbons in $\G$ on the  boundary components of $\G$ and the corresponding state of $\R_G$. 
Suppose that a ribbon $e$ of $\G$  intersects the discs of $\G$ in the arcs $a = [a_1,a_2 ] $ and $b=[b_1, b_2]$, and that $a_1$ and $b_1$ are connected along the ribbon $e$, and $a_2$ and $b_2$ are connected along the ribbon $e$. Also denote the corresponding points in $\s_{E(\G)}$ by $a_1, a_2, b_1$ and $b_2$.
There are two cases to consider.

Case 1: Removing  ribbons $e$ with $+$ weight from $\G$.   Before the removal, the boundary of $\G$ contains arcs  $[a_1 , b_1]$  and $[a_2,b_2]$, and after the removal of $e$, the boundary  contains arcs  $[a_1 , a_2]$  and $[b_1,b_2]$:
\[ 
\includegraphics[height=2cm]{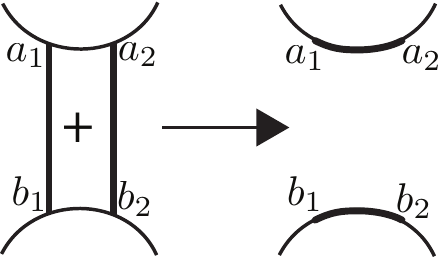}
 . \] 

Since $e$ is positive, the removal of an ribbon $e$ of $\G$ corresponds to the addition of a ribbon to $\s_{E(\G)}$. In $\s_{E(\G)}$, before the addition of this ribbon, $\s_{E(\G)}$ contains arcs  $[a_1 , b_1]$  and $[a_2,b_2]$, and after the addition of the corresponding ribbon, the boundary  contains arcs  $[a_1 , a_2]$  and $[b_1,b_2]$:
\[ 
\includegraphics[height=2cm]{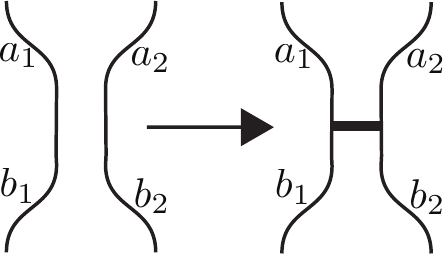}
 .\] 

 Case 2: Removing  ribbons $e$ with $-$ weight from $\G$.   Before the removal, the boundary of $\G$ contains arcs  $[a_1 , b_1]$  and $[a_2,b_2]$, and after the removal of $e$, the boundary  contains arcs  $[a_1 , a_2]$  and $[b_1,b_2]$:
\[ 
\includegraphics[height=2cm]{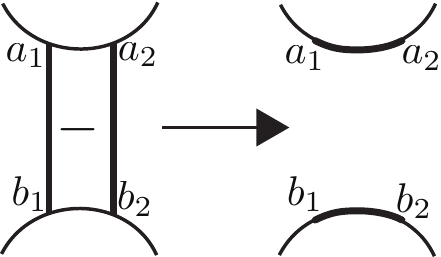}
. \]

Since $e$ is negative, the removal of an  ribbon $e$ of $\G$, corresponds to the removal of  an ribbon of $\s_{E(\G)}$. In $\s_{E(\G)}$, before the removal of the  ribbon, $\s_{E(\G)}$ contains arcs  $[a_1 , b_1]$  and $[a_2,b_2]$, and after the removal of the corresponding ribbon, the boundary  contains arcs  $[a_1 , a_2]$  and $[b_1,b_2]$:
\[ 
\includegraphics[height=2cm]{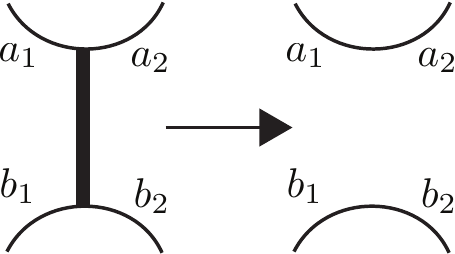}
.  \]

By combining these two observation on the effect of the addition and removal of ribbons on the number of boundary components of $\G$ and and $\s_{E(\G)}$, we see that when  $A\subset E(\G)$, $\partial(A) = \partial (\s_A)$. Thus for a fixed $A\subset E(\G)$,
\[   a^{\partial (A)} \prod_{e\in A} y_e = 
 (y_+)^{|E_+|}  a^{\partial (\s_A)} (y_-)^{|\s_A|} \]
and the result follows.
\end{proof}
The proof of the proposition then follows easily:
\begin{proof}[Proof of Proposition~\ref{pr:unsig}]
By Lemma~\ref{lem:pp1}, 
$ Z_{Potts}(G) =  q^{v(\G)/2}  \sum_{A \subset \E(\G)}  q^{\partial (A)/2}\prod_{e\in A}x_{e}$, which, by an application of Lemma~\ref{lem:pp2}, is equal to 
$  q^{|V|/2} x_+^{|E_+|} \sum_{ \A \subset E(\R_G)  } q^{\partial (A) /2} x_-^{|A|} $. 
\end{proof}

\subsection{The embedded graph $\M_G$}\lbl{ss:M}
In this subsection we construct a  graph $\M_G$ embedded on a surface $\Sigma_G$. Later we will construct an unsigned vertex model for $Z_{Potts}$ using this embedded graph. 

Consider the ribbon graph $\R_G$. As mentioned in Section~\ref{sec:pre}, 
we may regard a ribbon graph as a punctured surface (with additional structure). 
We define the closed, orientable surface $\Sigma_G$ to be the surface obtained from $\R_G$  by capping  the punctures with discs. 
Observe that there is a natural embedding $\R_G \subset \Sigma_G$. 

Construct an embedded graph   $\M_G \subset \R_G \subset  \Sigma_G$ as follows:
\begin{enumerate}
\item Place a vertex $v_e$ at the image of the mid-point $\{1/2 , 1/2\}$ of each ribbon $e$ of $\R_G$ under its natural embedding in $\Sigma_G$.
\item  Place vertices $v_{e_1},v_{e_2},v_{e_3},v_{e_4}$ at the images of the four points $\{0\}\times \{0\}, \{0\}\times \{1\}, \{1\}\times \{0\}, \{1\}\times \{1\}$ of each ribbon $e$ of $\R_G$ under its natural embedding in $\Sigma_G$. (These are the points where the corners of the  ribbons and discs of $\R_G$ intersect.)
\item Add edges $(v_e , v_{e_i})$, for $i=1,2,3,4$ and each ribbon  $e\in E(\R_G)$. We assume each edge is contained in the image of the ribbon $e$.
\item In each disc $v$ of $\R_G\subset \Sigma_G$, whenever two vertices $v_{e_i}$ and $v_{f_j}$ are connected by a path on  $\partial(v) - \{ e \in E(\R_G)   \}$, add an edge $(v_{e_i},v_{f_j})$. We assume this edge is contained in the disc $v$  
\end{enumerate}
We call the graph $\M_G$  the {\em unsigned medial graph of $\R_G$}.

An example of this construction is shown in Figure~\ref{fig:medexamp}, where $\R_G$   is embedded on the torus and is shown on the left and the medial graph is shown in the middle of the figure.
\begin{figure}
\[\begin{array}{ccc}
\includegraphics[height=3cm]{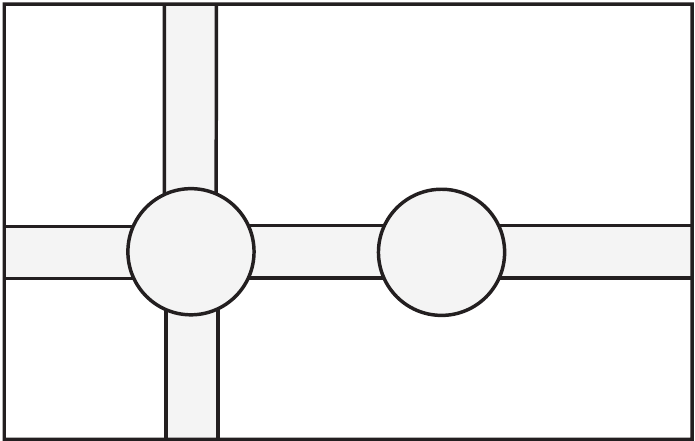}  &
\includegraphics[height=3cm]{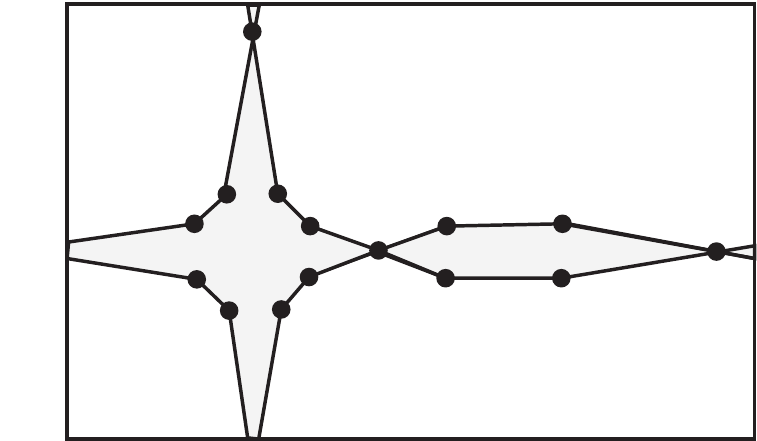}
&
\includegraphics[height=3cm]{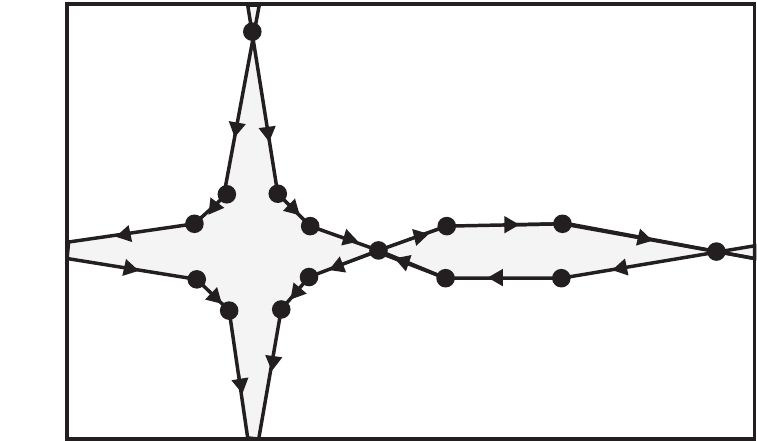}
\\
\R_G & \M_G & \text{An element of } AC(\M_G) 
  \end{array}\]
\caption{An example of the medial ribbon graph $\M_G$}
\lbl{fig:medexamp}
\end{figure}

Notice that each vertex of $\M_G$ is of valency two or four and that locally the set of   edges incident   with each $v_e$ meet transversally. 

We will give the embedded graph $\M_G$ the checkerboard colouring in the following way.
Colour the regions of $\M_G \subset \Sigma_G$ black or white according to the following rules: colour all ribbons and discs of the embedded ribbon graph $\R_G \subset \Sigma_G$ black. Since   $\M_G \subset \R_G \subset \Sigma_G$, some of the regions of $\M_G$ are coloured entirely black. By colouring all of the remaining regions  (i.e. those that are not coloured entirely black) white we obtain our colouring.

\medskip

 By an  {\em A-smoothing} of $\M_G$ at $v_e$ we mean the  local replacement of configurations:
\[ 
\includegraphics[height=1.5cm]{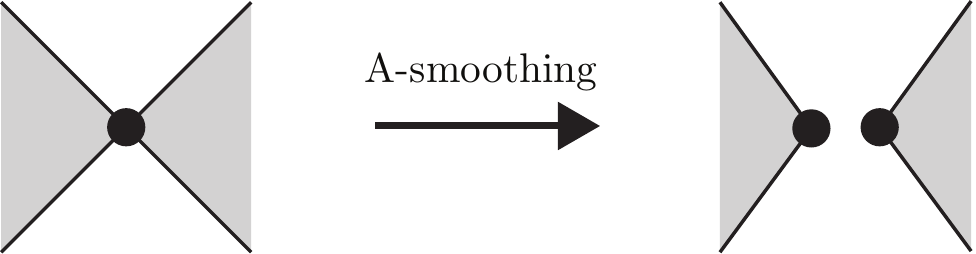}
 ;\]
and by a   {\em B-smoothing} of $\M_G$ at $v_e$ we mean the  local replacement of configurations:
\[
\includegraphics[height=1.5cm]{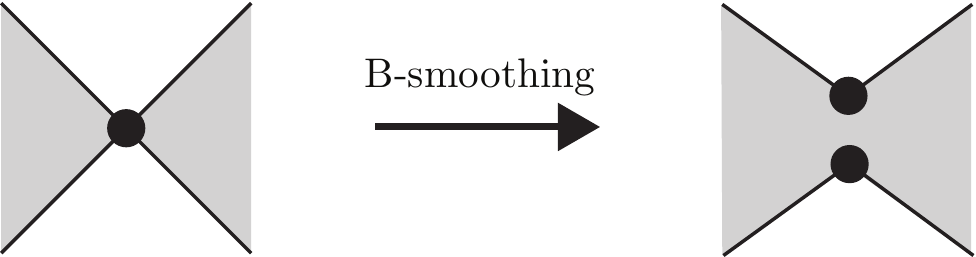}
.\]

Notice also that an A-smoothing splits black regions of the colouring, and a B-smoothing joins black regions of the colouring.

A {\em smoothing} of $\M_G$ is the graph obtained by a smoothing of each of the four-valent vertices  $v_e$, $e \in E(\R_G)$. A {\em cycle} in a smoothing  is a connected component of the smoothing. Therefore, the number of cycles of a smoothing of $\M_G$ is just the number of connected circles left over when all the vertices of $\M_G$ are smoothed.
We denote the set of smoothings of $\M_G$ by $S(\M_G)$.

\begin{lemma}\lbl{lem:smoo}
If    $x_+=(x_-)^{-1}$, where $x_{\pm}=q^{-1/2}w_{\pm}$, then 
\[  Z_{Potts} (G) =  q^{V/2} x_+^{|E_+|} \sum_{ s \in S(\M_G)  } q^{p(s) /2} x_-^{B(s)}, \]
where $p(s)$ is the number of cycles in the smoothing $s$, and  $B(s)$ denotes the number of B-smoothings used in the construction of $s$.
\end{lemma}
\begin{proof}
By proposition~\ref{pr:unsig}, it is enough to show that 
\[  \sum_{ \A \subset E(\R_G)  }  (q^{1/2})^{\partial (\A) } x_-^{|\A|}    =  \sum_{ s \in S(\M_G)  } (q^{1/2})^{p(s) } x_-^{B(s)} .\]
To see why this equality holds consider the following  bijection between the spanning subgraphs $(V(\R_G) , \A)$ of $\R_G$ and the smoothings of $\M_G$.
Let $E(\R_G)=\{ e_1 ,e_2 , \ldots , e_k\}$. Also let $\{ v_1 ,v_2 , \ldots , v_k\}$ be the vertices of  $\M_G$ such that $v_i$ arises at the center of $e_i$ in the construction of $\M_G$ from $\R_G$. 
Given $\A \subset E(\R_G)$, construct a smoothing $s$ of $\M_G$ by B-smoothing at vertices $v_i$ when $e_i \in \A$, and A-smoothing at vertices $v_i$ when $e_i \notin \A$. Clearly, $|\A|= B(s)$. 

Let $s$ be a smoothing  and $\A$ be the corresponding spanning ribbon sub-graph. We need to show that $\partial(\A) = p(s)$. We recall that $\M_G\subset \R_G$ and therefore $s\subset \R_G$. Deform the embedded graph $s$ by ``sliding'' its edges and vertices along the edges of $\R_G$ so that $s$ is embedded in the set of boundaries of $V(\R_G) \cup E(\R_G)$. It is readily seen that $s$ lies on the boundary $I \times \{0,1\}$ of an edge $e$ if and only if $e\in \A$; and that $s$ lies  on the boundary $ \{0,1\}\times I \subset V(\R_G)$ of an edge $e$ if and only if $e\notin \A$. It then follows that $s$ embeds in the set of boundary components of the spanning ribbon sub-graph $(V(\R_G),\A)$ and that $s$ necessarily covers the boundary components of $(V(\R_G),\A)$. Thus  $\partial(\A) = p(s)$, and the equality of the state sums then follows.
\end{proof}

\section{The Kauffman bracket}\lbl{sec:kauffman}

\subsection{The Kauffman bracket as a Potts model}

Let $D\subset \mathbb{R}^2$ be a link diagram. 
By the {\em checkerboard colouring} of $D$, we mean
the  assignment of  a colour black or white to each region of $D$ in such a way  that adjacent regions have different colours and the unbounded region is white. The {\em Tait sign}  of a crossing of $D$ is an element of $\{+,-\}$ which is assigned to the crossing according to the following scheme:
\[ \begin{array}{ccc}
\raisebox{1mm}{\includegraphics[height=1.5cm]{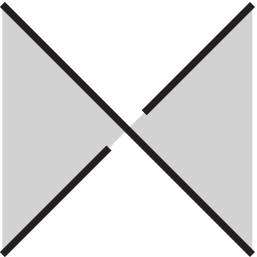}}
& \quad& 
\raisebox{1mm}{\includegraphics[height=1.5cm]{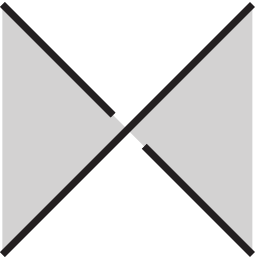}}
\\
+ & & -
\end{array}.\]

 The {\em Tait  graph} $T_D$ is a signed plane graph constructed from $D$ as follows: given the link diagram,  place a vertex in each black region, and add an edge between two vertices whenever the corresponding regions of $D$ meet at a crossing. Colour each edge of the graph by the medial sign of the corresponding crossing.

It follows from \cite{Ka3} that the Kauffman bracket of a link diagram $D$ can be calculated from its Tait graph $T_D=(V,E)$ by the Potts partition function:
\begin{equation}\label{eq:kp}   \langle D \rangle = \delta^{|V|-1} A^{|E_-|-|E_+|} \sum_{A \subset E} (\delta^2)^{k(A)} \prod_{e\in A} \delta A^{2b_e},  \end{equation}
where $\delta = -A^2-A^{-2}$ and $b_e = \pm 1$ if $w(e)=\pm$. (See also  \cite{Jo, Lo, Wu1} for equivalent expressions of the Jones polynomial and Kauffman bracket as a Potts model and \cite{Wu2} for an alternative expression.)

\subsection{Connections with links in thickened surfaces}

Let $\fL$ be a link in a thickened surface $\Sigma\times I$. The {\em Kauffman bracket} $\langle \fD \rangle$ of a diagram $\fD$ of $\fL$ is defined (see \cite{ik}) as  the state sum 
\begin{equation} \label{eq:kbthick} 
\langle \fD \rangle = \sum_{\fs \in S(\fD)}  A^{\alpha(\fs)} A^{- \beta (\fs)} \delta^{c(\fs)-1},
\end{equation} 
where a state is an A-smoothing or B-smoothing of each crossing: 
\[
\includegraphics[height=1.5cm]{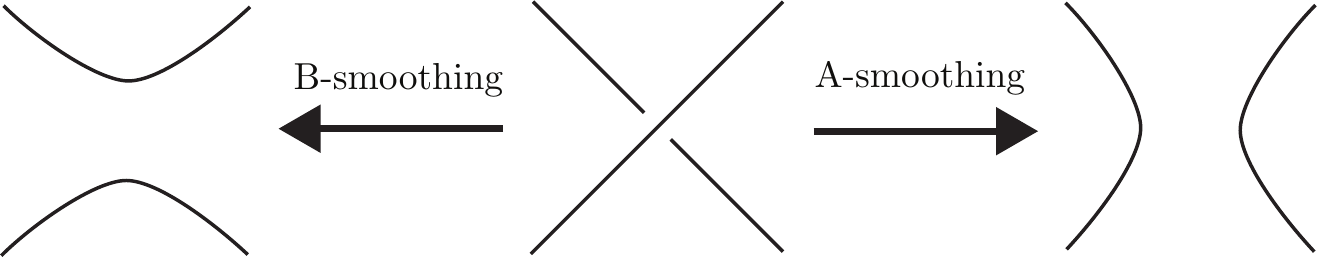}
;\]
the sum is over all states $S(\fD)$ of $\fD$, $\alpha(\fs)$ is the number of A-smoothings used to construct the state; $\beta(\fs)$ is the number of B-smoothings used to construct the state;  $c(\fs)$ is the number of cycles in a state; and $\delta = -A^2-A^{-2}$.

The following theorem says that for every link $L$ in $S^3$ with diagram $D$, there exists an alternating link $\fL$ in a thickened surface with diagram $\fD$ such that  $\langle D \rangle =\langle \fD \rangle$.
\begin{theorem}\lbl{th:twojones}
Let $D$ be a link diagram and $\M_{D}\subset \Sigma \times I$ be the unsigning of the medial graph $M_D$ of $D$. Let $\fD_D \subset \Sigma$  be the link diagram obtained from $\M_D$ by replacing each four-valent vertex with a crossing of  negative medial sign.
Then 
\[  \langle D \rangle =\langle \fD_D \rangle    .\]
\end{theorem} 
\begin{proof}
First consider the 
Kauffman bracket $ \langle D \rangle $. If   $M_D= (V,E)$.
\[   \langle D \rangle = \delta^{|V|-1} A^{|E_-|-|E_+|} \sum_{\A \subset E} (\delta^2)^{k(\A)} \prod_{e\in \A} \delta A^{2b_e},  \]
where $\delta = -A^2-A^{-2}$ and $b_e = \pm 1$ if $w(e)=\pm$. 
Applying Proposition \ref{pr:unsig} to this expression gives
\[   \langle D \rangle  =  
\delta^{|V|-1} A^{|E_-|-|E_+|}
 \delta^{V} (A^{2} )^{|E_+|} \sum_{ \A \subset E(\R_{M_G})  } \delta^{\partial (\A)} A^{-2|\A|}
\]
or 
\begin{equation}\label{eq:tj2}    \langle D \rangle  =  
\delta^{-1} A^{|E|}
   \sum_{ \A \subset E(\R_{M_G})  } \delta^{\partial (\A)}  A^{-2|\A|}
\end{equation}

Next consider the unsigned ribbon graph $\R$ which is the unsigning of  the signed medial graph $M_D$ of $D$. Construct a checker-board coloured link diagram $\fD$ from $\R$ by replacing each ribbon with a crossing of negative medial sign; connecting the crossings according to the boundary of $\R$; and colouring regions of $\fD$ which contained the interiors of the discs of $\R$ black.  Clearly  $\fD$ and the diagram $\fD_D$ from the statement of the theorem are related by isotopy of $\Sigma$. 

Every disc of $\R$ corresponds with a black coloured region of $\fD$. Moreover  an A-smoothing of $\fD$ splits  black regions, B-smoothing of $\fD$ joins black regions. 
If we construct a spanning ribbon sub-graph of  $\R$   from a smoothing of $\fD$ by including a ribbon if and only if the corresponding crossing is B-smoothed, then the number of boundary components of this spanning ribbon sub-graph is equal to the number of cycles in the smoothing and the number of edges corresponds with  the number of B-smoothings used in the smoothing of $\fD$. Therefore 
\[\langle \fD \rangle = \sum_{\fs}  A^{\alpha(\fs)} A^{- \beta (\fs)} \delta^{c(\fs)-1}
= \sum_{\A \subset E(\R)} A^{|E (\R) |-|\A|} A^{-|\A|}\delta^{\partial (\A)-1}
\]
Thus 
\begin{equation}\label{eq:tj1} \langle \fD_D \rangle  = \delta^{-1} A^{|E(\R)|}
   \sum_{ \A \subset E(\R)  } \delta^{\partial (\A)}  A^{-2|\A|}. \end{equation}

 By comparing (\ref{eq:tj2}) and (\ref{eq:tj1}) we see
 $ \langle D \rangle =\langle \fD_D \rangle  $ as required.

\end{proof}

It is a straightforward exercise to rephrase this theorem in terms of virtual links rather than links in a thickened surface (c.f. \cite{CP1} and \cite{CP}). 

\subsection{The connection with the \BR polynomial}\lbl{ss:graphpolys}
Part of our motivation for this work came from expressions for the Jones polynomial as evaluations of    the Tutte polynomial of a graph and
the \BR polynomial \cite{BR1, BR} of a ribbon (or embedded) graph. Here we recover and relate two results from the literature which express the Jones polynomial in terms of the \BR polynomial.

The {\em \BR polynomial for ribbon graphs} is defined as the sum
\[
R(\F ;X,Y,Z)  = \sum_{\A \subset E(\F)}(X-1)^{r(\F)-r(\A)}Y^{n(\A)}Z^{k(\A)-\partial(\A)+n(\A)}.
\]
Observe that when  $Z =1$ this polynomial becomes the Tutte
polynomial:
\[
R(\F; X ,Y-1 ,1) = T(\F;X,Y) :=  \sum_{\A \subset E(\F)}(X-1)^{r(\F)-r(\A)}(Y-1)^{n(\F)} .
\]
This coincidence of polynomials also holds when the ribbon graph $\F$
is planar (since the largest  exponent of $Z$ is twice the genus the surface $\F$).

\medskip

A key observation in the proof of Theorem~\ref{th:twojones} was that we can regard the ribbon graph $\R_D$ as either the unsigned ribbon graph of a link diagram $D$, or as the Tait graph of a link diagram in a thickened surface. By taking each of these two perspectives in turn, we can recover and relate the results from \cite{Da} and \cite{CP1} which relate the \BR polynomial and the Kauffman bracket.

The \BR polynomial is given by 
\[R(\R;X,Y,Z)= \sum_{\A \subset E(\R)} (X-1)^{r(\R)-r(\A)}Y^{n(\A)}Z^{k(\A)-\partial (\A)+n(\A)}.\] 
We may rewrite this as
\begin{equation}\label{eq:expbr} R(\R;X,Y,Z)= (X-1)^{-k(\R)} (YZ)^{-|V(\R)|}  
\sum_{\A \subset E(\R)} [(X-1)YZ^2]^{k(\A)}(YZ)^{|\A|}Z^{-\partial (\A)}.  \end{equation}
Now consider the state sum 
\[  \delta^{-1} A^{|E|}
   \sum_{ \A \subset E(\R_{M_G})  } \delta^{\partial (\A)}  A^{-2|\A|} .\] 
Setting $X=-A^4$, $Y=A^{-2}\delta$ and $Z=\delta^{-1}$ in (\ref{eq:expbr}) gives
\[   \delta^{-1} A^{|E|}
   \sum_{ \A \subset E(\R_{M_G})  } \delta^{\partial (\A)}  A^{-2|\A|}  =  
\delta^{-1+k(\R)} A^{|E|+2k(\R)-2|V(\R)|}
R(\R_{M_G} ; -A^4 , A^{-2}\delta , \delta^{-1}).
\]
$\R$ is connected so this expression simplifies to give
\begin{equation}\label{eq:jbr1}   \delta^{-1} A^{|E|}
   \sum_{ \A \subset E(\R_{M_G})  } \delta^{\partial (\A)}  A^{-2|\A|}  =  
 A^{|E|+2-2|V(\R)|}
R(\R_{M_G} ; -A^4 , A^{-2}\delta , \delta^{-1}).
\end{equation}

Viewing $\R_D$ as  the unsigning of the signed  graph of a link diagram $D$, and substituting (\ref{eq:jbr1}) into Equation~(\ref{eq:tj1}), recovers:
\begin{theorem}[\cite{Da}]
Let $D$ be a link diagram and $\R$ be the unsigning of the signed  graph of a link diagram $D$, then
\[ \langle D \rangle =   
 A^{|E|+2-2|V(\R)|}
R(\R ; -A^4 , A^{-2}\delta , \delta^{-1}).\]
\end{theorem}
Viewing 
 $\R_D$ as the  Tait graph of a link diagram in a thickened surface, and substituting (\ref{eq:jbr1}) into Equation~(\ref{eq:tj2}), recovers
\begin{theorem}[\cite{CP}]
Let $\R$ be the Tait graph of an alternating link diagram $D$ on a surface, then 
\[ \langle D \rangle =  
 A^{|E|+2-2|V(\R)|}
R(\R ; -A^4 , A^{-2}\delta , \delta^{-1}).\]
\end{theorem} 
We see that these two theorems from \cite{CP} and \cite{Da} which relate the \BR polynomial and the Kauffman bracket follow from two different interpretations of the  same ribbon graph.


 \section{An unsigned vertex model for the Jones polynomial}
In this section we show that the Jones polynomial has a construction as a vertex model of an unsigned graph. The unsigned graph is  the checkerboard coloured unsigned medial graph $\M_G$ which is embedded in some surface $\Sigma_G$.

\subsection{Arrow coverings}\lbl{ss:arrow}
\begin{definition}
An {\em arrow covering} of $\M_G$ or a smoothing of $\M_G$ is an assignment of an orientation to each edge such that every vertex has an equal number of incoming and outgoing edges.

The set of arrow coverings of $\M_G$ or of a smoothing $s \in S(\M_G)$ will be denoted by
$AC(\M_G)$ or  $AC(s)$ respectively. 
\end{definition}
An example of an arrow covering of $\M_G$ is shown on the right of Figure~\ref{fig:medexamp}.

Our aim is to express  $Z_{potts}(G) $ as a state sum over arrow coverings of $\M_G$, where the terms in this state sum arise as weights associated with the vertices of $\M_G$.

We  take Lemma~\ref{lem:smoo}, which states that 
\[   Z_{Potts} (G) =  q^{V/2} x_+^{|E_+|} \sum_{ s  \in S(\M_G) } q^{p(s) /2} x_-^{B(s)},\]
 as our starting point and rewrite this state sum as a vertex model.
 
Our first step towards  this goal  is to  rewrite the expression $q^{p(s) /2}$, where $s$ is a smoothing of  $\M_G$ and $p(s)$ is the number of cycles in $s$, in terms of arrow coverings of $s$.
To do this we define a {\em cycle weight function}, $\delta$, to be a function from closed oriented curves on a surface $\Sigma$ to the set $\{+1,-1\}$ which satisfies the following condition: if $c$ is an oriented curve on $\Sigma$ and $-c$ is the curve with the opposite orientation, then $-\delta(c)=\delta(-c)$. 
We say that a curve $c$ is {\em positively weighted} if $\delta(c)=+1$, and is {\em negatively weighted} if $\delta(c)=-1$.

Our vertex model uses a choice of cycle weight function, however the construction is independent of the actual choice of the cycle weight function. For completeness we define a  cycle weight function.
Let $c\subset \Sigma_G$ be a closed oriented curve. Then either $c$ bounds a disc in $\Sigma_G$, or $c$ is homologically non-trivial, or $c$ separates the surface $\Sigma_G$ into two components. The surface  $\Sigma_G$  is either a 2-sphere or the connected sum of $g$ tori. If the surface is the connected sum of tori then label the tori $1,2,\ldots, g$.  
 Also choose an ordered set of generators $g_1, \ldots , g_{k}$ of $H_1 (\Sigma_G)$. A cycle weight function can be defined in the following way.
First, if $c$ bounds a disc in $\Sigma_G$ then take $\delta(c)$ to be the rotation number of the tangent vector  of $c$ in the tangent space as we travel around the curve $c$ in the direction of its orientation. Secondly, if $c$ is homologically non-trivial  then we can write the homology class  $[c]$ of $c$ in terms of our generators: $[c]=a_1g_1+\cdots +a_kg_k$. We then define $\delta (c)$ to be the sign of the first non-zero coefficient $a_i$ in this expression.
Finally, if $c$ is null-homologous and does not bound a disc in $\Sigma_G$, then cut the surface $\Sigma_G$ along $c$ to obtain two punctured surfaces. Cap-off the punctures using discs to obtain two closed surfaces each containing a copy of $c$. Let $\delta(c)$ be the winding number of the copy of the curve $c$ contained in the surface that contains the torus labeled $1$.
 It is clear that $\delta(-c)=-\delta(c)$ in all three situations and so $\delta$ is a cycle weight function.    

\begin{lemma}\lbl{lem:cwf}
Let $\delta$ be a cycle weight function and $s$ be a smoothing of $\M_G$. Then 
\[  \sum_{s\in S(\M_G)} r^{B(s)}  (t+t^{-1})^{p(s)} = \sum_ {s\in S(\M_G)}   r^{B(s)}  \sum_{a\in AC(s)}  \prod_{c\in a}  t^{\delta(c)},  \]
where $AC(s)$ denotes the set of arrow coverings of a smoothing $s$ and the product is over all oriented cycles in $a$.
\end{lemma}
\begin{proof}
To prove the lemma we need to show \[   \sum_{a\in AC(s)}  \prod_{c\in a} t^{\delta(c)}  = (t+t^{-1})^{p(s)}. \]

Consider the expression $ \sum_{a\in AC(s)}  \prod_{c\in a}  t^{\delta(c)}$. 
Suppose that there are $c$ cycles in the state $s$. Then an arrow covering $a\in AC(s)$ orients $k$ of the cycles negatively and $c-k$ cycles positively. Such an arrow covering contributes the summand $t^{c-k}t^{-k}$ to $ \sum_{a\in AC(s)}   \prod_{c\in a} t^{\delta(c)}$.

Since the orientations in an arrow covering are independent of each other, for each $k=0 , \ldots ,c$ there are exaclty $\binom{c}{k}$ choices for the $k$ negatively weighted cycles. Therefore 
\[  \sum_{a\in AC(s)}   \prod_{c\in a} t^{\delta(c)}= \sum_{k=0}^c t^{c-k}t^{-k} =  (t+t^{-1})^{c}.\]
 The result then follows since $c=p(s)$.
\end{proof}

The next step in the construction of our vertex model involves writing the sum from Lemma~\ref{lem:cwf} as a sum over arrow coverings.

\medskip

There is a natural map 
\[  \xi :  \{  a \; | \; a\in AC(s), \;\; s\in   S(\M_G) \}    \longrightarrow  AC(\M_G)  \]
from the set of  all arrow coverings of all smoothings of $\M_G$ to the set of all arrow coverings of $\M_G$, given by
\[    
\includegraphics[height=1cm]{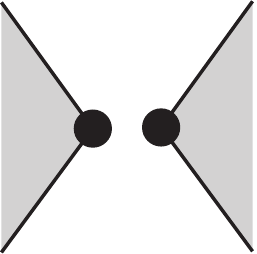}
  \raisebox{5mm}{\;\;$\mapsto$\;\;}  
  \includegraphics[height=1cm]{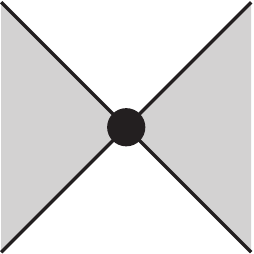}
  \hspace{1cm}   \raisebox{5mm}{\text{and}}  \hspace{1cm}  
\includegraphics[height=1cm]{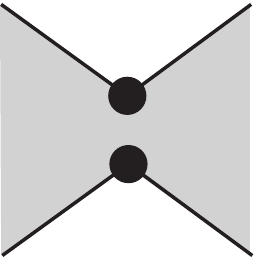}
  \raisebox{5mm}{\;\;$\mapsto$\;\;} 
  \includegraphics[height=1cm]{figs/xi4}
   \]
 Where the orientations of the edges (not shown in the figure) are unchanged.

This map is easily seen to be surjective and, in general, not injective.  
We  define a map $\rho$ from $AC(\M_G)$ to $\mathbb{Z}[r]$-linear combinations of arrow coverings of smoothings of $\M_G$ by the figure:
\[ \begin{array}{ccl}
\includegraphics[height=1cm]{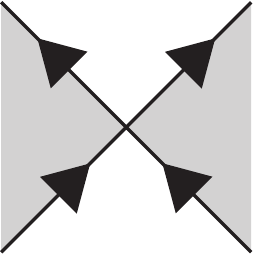}
 &  \raisebox{5mm}{$\overset{\rho}{\mapsto}$}&
 \includegraphics[height=1cm]{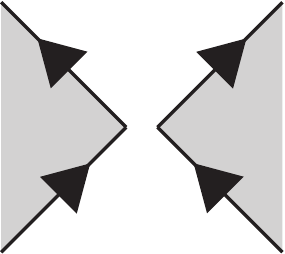}
 \\
\includegraphics[height=1cm]{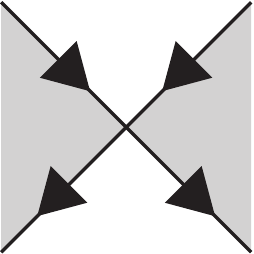}
 &  \raisebox{5mm}{$\overset{\rho}{\mapsto}$}&
  \includegraphics[height=1cm]{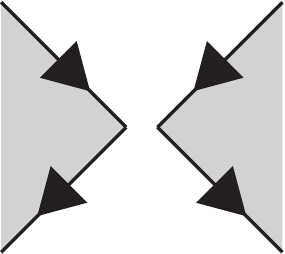}
  \\
\includegraphics[height=1cm]{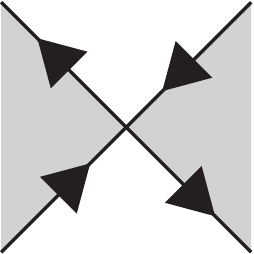}
&  \raisebox{5mm}{$\overset{\rho}{\mapsto}$}&
  \includegraphics[height=1cm]{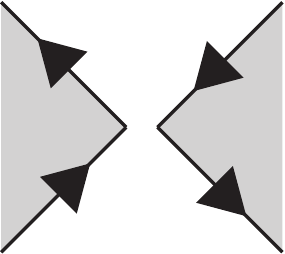}
    \raisebox{5mm}{$\;+r\;$} 
    \includegraphics[height=1cm]{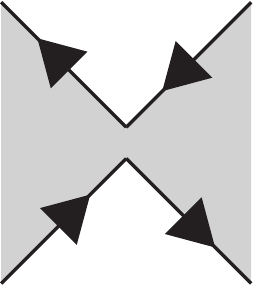}
\end{array} 
 \hspace{2cm}
 \begin{array}{ccl}
\includegraphics[height=1cm]{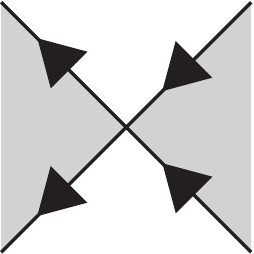}
&  \raisebox{5mm}{$\overset{\rho}{\mapsto}$}&
 \raisebox{5mm}{$r\;$} 
 \includegraphics[height=1cm]{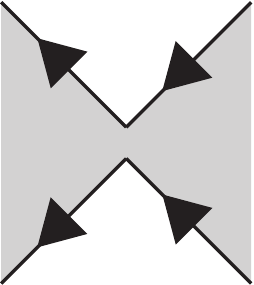}
   \\
\includegraphics[height=1cm]{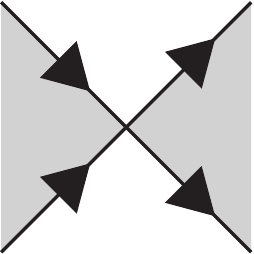}
&  \raisebox{5mm}{$\overset{\rho}{\mapsto}$}&
 \raisebox{5mm}{$r \;$} 
 \includegraphics[height=1cm]{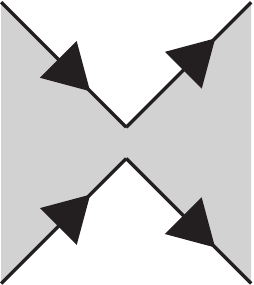}
  \\
\includegraphics[height=1cm]{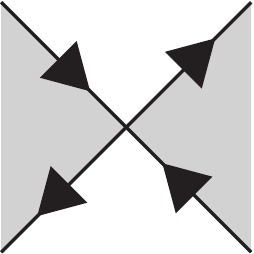}
&  \raisebox{5mm}{$\overset{\rho}{\mapsto}$}&
  \includegraphics[height=1cm]{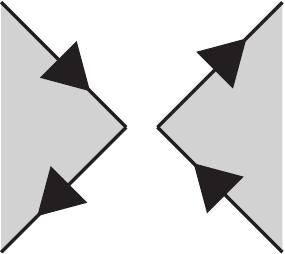}
   \raisebox{5mm}{$\;+r\;$}
   \includegraphics[height=1cm]{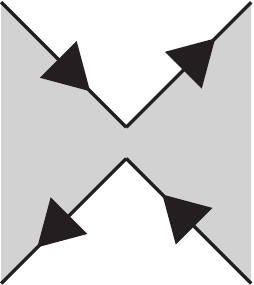}
\end{array} . \]
(For clarity, the vertices are not shown in the figure.)

In the set $\rho( AC(\M_G) )$ 
every  element of $ \{  a   \; | \; a\in AC(s), \;\;s\in   S(\M_G) \}$  occurs exactly once as a summand (to see why this is consider the inverse image of the map $\xi$). Also
  every summand  $ac(s)$ in  $\rho( AC(\M_G) )$ is weighted  by $r^{B(s)}$.
 These two observations form the basis of the proof of the following lemma.
 \begin{lemma}\lbl{lem:rho}
 \[  \sum_ {s\in S(\M_G)}   r^{B(s)}  \sum_{a\in AC(s)}  \prod_{c\in a}  t^{\delta(c)} = \sum_{a\in AC(\M_G) } W(\rho(a)) ,  \]
 where  $W$ is the function from $\mathbb{Z}[r]$-linear combinations of oriented cycles to $\mathbb{Z}[r,t,t^{-1} ]$, which is defined as the linear extension of 
 \[ W:   \bigcup_{i} c_i  \mapsto  \prod_i t^{\delta(c_i)}.  \]
 \end{lemma}
  \begin{proof}
A summand in the expression $ \sum_ {s\in S(\M_G)}   r^{B(s)}  \sum_{a\in AC(s)}  a$ consists of an arrow covering of a state $s$ multiplied by $r^{B(s)}$. Moreover every arrow covering of every state arises in exactly one summand.

A summand in $\sum_{a\in AC(\M_G) } \rho(a)$ also consists of an arrow covering of a state $s$ multiplied by $r^{B(s)}$. Again  every arrow covering of every state arises in exactly one summand. 

Therefore, 
\[  \sum_ {s\in S(\M_G)}   r^{B(s)}  \sum_{a\in AC(s)}  a = \sum_{a\in AC(\M_G) } \rho(a),  \] 
Since the terms $  \prod_{c\in a}  t^{\delta(c)}$ on the left-hand side  and $W(\rho(a))$ on the right-hand side both map a collection of oriented cycles $ \bigcup_{i} c_i$ to $ \prod_i t^{\delta(c_i)}$
 \[  \sum_ {s\in S(\M_G)}   r^{B(s)}  \sum_{a\in AC(s)}  \prod_{c\in a}  t^{\delta(c)} = \sum_{a\in AC(\M_G) } W(\rho(a)) .  \]
  \end{proof}

We are now in a position to state and prove our expression for  the Potts partition function of a signed plane graph  as an unsigned vertex model.
\begin{theorem}\lbl{th:main}
If    $x_+=(x_-)^{-1}$, where $x_{\pm}=q^{-1/2}w_{\pm}$, and $q={t+t^{-1}}$ , then 
\[  Z_{Potts} (G) =  (t+t^{-1})^{V/2} x_+^{|E_+|} \sum_{a\in AC(\M_G) } W(\rho(a)) , \]
where  $W$ is the function from $\mathbb{Z}[x_-]$-linear combinations of oriented cycles to $\mathbb{Z}[x_-,t,t^{-1} ]$ which is defined as the linear extension of 
 \[ W:   \bigcup_{i} c_i  \mapsto  \prod_i t^{\delta(c_i)}.  \]
\end{theorem}
\begin{proof}
The theorem follows easily  from Lemmas~\ref{lem:smoo}, \ref{lem:cwf} and \ref{lem:rho}.
\end{proof}

\medskip

We know how to express the Kauffman bracket of a link diagram as a Potts model (see Section~\ref{sec:kauffman}). We may apply Theorem~\ref{th:main} to this  expression for the Kauffman bracket and  write the Kauffman bracket of any planar link diagram as a vertex model of an unsigned graph.

\begin{theorem}\lbl{th:kauff1}
Let $D$ be a link diagram on $\mathbb{R}^2$ and $T$ be its Tait graph. Then if $\M$ is the medial graph of the unsigning $\R$ of $T$, the Kauffman bracket can be recovered from the vertex model of an unsigned embedded graph $\M$
\[ \langle D \rangle =  (-A^2-A^{-2})^{-1} A^N   \sum_{a\in AC(\M) } W(\rho(a)) ,
   \]
where $N$ is the number of crossings of $D$ and 
$W$ is the function from $\mathbb{Z}[A^{-2}]$-linear combinations of oriented cycles to $\mathbb{Z}[A, A^{-1} ]$ which is defined as the linear extension of 
 \[ W:   \bigcup_{i} c_i  \mapsto  \prod_i (-A^2)^{\delta(c_i)}.  \]

\end{theorem}
\begin{proof}
 Equation~(\ref{eq:tj2}) gives
\[  \langle D \rangle = (-A^2-A^{-2})^{-1} A^{|E|}
   \sum_{ \A \subset E(\R)  }(-A^2-A^{-2})^{\partial (\A)}  A^{-2|\A|}. \]
By Proposition~\ref{pr:unsig} and Lemma~\ref{lem:smoo}, we can write this as
\[  \langle D \rangle = (-A^2-A^{-2})^{-1} A^{|E|}
   \sum_{ s \in S (\M)  }(-A^2-A^{-2})^{p(s)}  (A^{-2})^{B(s)}. \]
By Theorem~\ref{th:main}, this is equal to the vertex model
\[  \langle D \rangle = (-A^2-A^{-2})^{-1} A^{|E|} \sum_{a\in AC(\M) } W(\rho(a))
  , \]
where $t=-A^2$ and $x_- = -A^{-2}$. The result then follows upon noting that $|E|=N$.

\end{proof}

We may also regard the  unsigned  graph $\M$ as the underlying  graph of an alternating link diagram on a surface. With this interpretation of $\M$, Theorem~\ref{th:kauff1} becomes: 
  \begin{theorem}\lbl{th:kauff2}
  Let $\M$ be the underlying graph of an alternating link diagram $\fD$ on a surface (with each crossing of negative medial sign).
 Then the Kauffman bracket $\langle \fD \rangle$ is given by the vertex model based on the graph $\M$: 
\[ \langle \fD \rangle =   (-A^2-A^{-2})^{-1} A^N   \sum_{a\in AC(\M) } W(\rho(a)) ,
   \]
where $N$ is the number of crossings of $\fD$ and 
$W$ is the function from $\mathbb{Z}[A^{-2}]$-linear combinations of oriented cycles to $\mathbb{Z}[A, A^{-1} ]$ which is defined as the linear extension of 
 \[ W:   \bigcup_{i} c_i  \mapsto  \prod_i (-A^2)^{\delta(c_i)}.  \]
\end{theorem}

The proof is similar to the proof of Theorem~\ref{th:kauff1}, except we start with Equation~(\ref{eq:tj1}), rather than Equation~(\ref{eq:tj2}). Due to this  similarity, the proof is omitted.



\begin{thebibliography}{99}


\bibitem{Ba} R. J. Baxter,  Exactly solved models in statistical mechanics. Academic Press, Inc., London, 1982.


\bibitem{BR1} B. Bollob\'{a}s and O. Riordan,  A polynomial for
graphs on orientable surfaces,  Proc. London Math. Soc.  \textbf{83}
(2001),  513-531.


\bibitem{BR} B. Bollob\'{a}s and O. Riordan,   A polynomial of graphs
on surfaces,  Math. Ann.  \textbf{323}  (2002),  no. 1, 81-96.


\bibitem{CP1} S. Chmutov and I. Pak, The Kauffman bracket and the Bollobas-Riordan polynomial of ribbon graphs, preprint,  {\tt arXiv:math.GT/0404475}.





\bibitem{CP} S. Chmutov and I. Pak,   The Kauffman bracket of virtual links and the Bollob‡s-Riordan polynomial, Mosc. Math. J. {\bf 7} (2007) 409-418, {\tt
arXiv:math.GT/0609012}.

\bibitem{CV} S. Chmutov, J. Voltz, Thistlethwaite's theorem for virtual links, J. Knot Theory Ramifications {\bf 17} (2008), no. 10, 1189-1198. {\tt arXiv:0704.1310 }.

\bibitem{Da} O. T. Dasbach, D. Futer, E. Kalfagianni, X.-S. Lin and  N. W. Stoltzfus,  The Jones polynomial and graphs on surfaces, J. Combin. Theory Series B, {\bf 98} (2) (2008), 384-399,  {\tt arXiv:math.GT/0605571}.

\bibitem{FK} Fortuin, C. M.; Kasteleyn, P. W. On the random-cluster model. I. Introduction and relation to other models.  Physica  {\bf 57}  (1972), 536-564.

\bibitem{HM} S. Huggett and I. Moffatt, Expansions for the \BR polynomial of separable ribbon graphs, preprint {\tt  arXiv:0710.4266}.

\bibitem{ik} K. Inoue and T. Kaneto,
A Jones type invariant of links in the product space of a surface and the real line.
J. Knot Theory Ramifications {\bf 3} (1994), no. 2, 153-161.


\bibitem{Jo2} V. F. R. Jones, 
 A polynomial invariant for knots via von Neumann algebras. Bull. Amer. Math. Soc. (N.S.) {\bf 12} (1985), no. 1, 103-111.

\bibitem{Jo} V. F. R. Jones, 
On knot invariants related to some statistical mechanical models. 
Pacific J. Math. {\bf 137} (1989), no. 2, 311-334. 

\bibitem{Ka} L. H. Kauffman, State models and the Jones polynomial. Topology {\bf 26} (1987), no. 3, 395-407.

\bibitem{Ka2} L. H. Kauffman, Statistical mechanics and the Jones polynomial. Braids (Santa Cruz, CA, 1986), 263-297, 
Contemp. Math., {\bf 78}, Amer. Math. Soc., Providence, RI, 1988. 


\bibitem{Ka3} L. H. Kauffman, 
 A Tutte polynomial for signed graphs,
Combinatorics and complexity (Chicago, IL, 1987). 
Discrete Appl. Math. {\bf 25} (1989), no. 1-2, 105-127. 



\bibitem{LM} M. Loebl and I. Moffatt, The chromatic polynomial of fatgraphs and its categorification, Advances in Mathematics, {\bf 217} (2008) 1558-1587,   {\tt arXiv:math.CO/0511557}.

\bibitem{Lo} M. Loebl, Chromatic Polynomial, q-Binomial Counting and Colored Jones Function, Advances in Mathematics, {\bf 211}, 546-565, {\tt arXiv:math/0412460}. 

\bibitem{Mo} I. Moffatt,  Knot invariants and the
Bollob\'{a}s-Riordan polynomial of embedded graphs, European J. Combin., {\bf 29} (2008) 95-107, {\tt arXiv:math/0605466}.


\bibitem{PeWu} J. H. H. Perk and F. Y. Wu, Graphical approach to the nonintersecting string model: star-triangle equation, inversion relation, and exact solution. Phys. A {\bf 138} (1986), no. 1-2, 100--124.

\bibitem{So} A. D. Sokal, The multivariate Tutte polynomial (alias Potts model) for graphs and matroids. Surveys in combinatorics 2005, 173-226, 
London Math. Soc. Lecture Note Ser., {\bf 327}, Cambridge Univ. Press, Cambridge, 2005, {\tt arXiv:math/0503607}.

\bibitem{Th} M. B. Thistlethwaite, 
A spanning tree expansion of the Jones polynomial, 
Topology {\bf 26} (1987), no. 3, 297-309. 


\bibitem{Tu} V. G. Turaev,  The Yang-Baxter equation and invariants of links, Invent. Math. {\bf 92} (1988), no. 3, 527-553.



\bibitem{Wu1}  F. Y. Wu,  Knot theory and statistical mechanics. Rev. Modern Phys. {\bf 64} (1992), no. 4, 1099-1131. 

 \bibitem{Wu2}  F. Y. Wu, Jones polynomial as a Potts model partition function. J. Knot Theory Ramifications {\bf 1} (1992), no. 1, 47--57.


\end{thebibliography}
\end{document}